\documentclass[12pt]{article}
\usepackage{amsmath}
\usepackage{amssymb}
\usepackage{geometry}
\usepackage{amscd}
\usepackage{xypic}
\usepackage{mathrsfs}
\usepackage{amsthm}

\newtheorem{lemma}{Lemma}[section]
\newtheorem{theorem}{Theorem}[section]
\newtheorem{definition}{Definition}[section]
\newtheorem{proposition}{Proposition}[section]

\theoremstyle{definition}
\newtheorem{example}{Example}[section]

\begin{document}

\title{Computing Moment Maps of Hypersurfaces using MAXIMA}
\author{Dun Liang}
\maketitle

\section{Introduction}
Without special announcements, all varietes and manifolds are defined over the field of complex numbers ${\mathbb C}$.

Geometric invariant theory (GIT see \cite{MFK}) answers the question whether quotients exist the category of algebraic varieties. The answer is that after deleting so called the unstable points for a specific linearization, one constructs a categorical quotient of a variety by the action of a reductive group and furthur more one gets a geometric quotient after selecting the stable points.

This paper discusses stability of the easiest non-trivial situation when $G$ acts linearly on the projective space $X={\mathbb P}^n({\mathbb C})$. Over the complex numbers, the problem can also be considered in a simplectic point of view when on regards $G$ as a complex Lie group and $X$ as a simplectic manifold, in this case the maximal compact subgroup $K$ of $G$ acts simplectically on $X$. Then there exists a moment map ${\frak m}:X\rightarrow {\rm Lie}(K)^\vee$. For each $x\in X$, fix an ${\rm SU}(V)$-invariant norm, the square length of the moment matrix $\|{\frak m}(x)\|^2$ is a Morse function, it determines the stability of $x$. 

MAXIMA is an open source computer algebra system, it has a language that is mature and friendly to mathematicians. But for algebro-geometric usage, not all functions are prepared in the library including some very elementary and important functions.

This paper constructs the mometent matrices of cubic curves and computes the cricical values of the Morse function $\|{\frak m}(x)\|^2$ as an example.
 Then general case of hypersurfaces can be computed simply by changing the announcement of the degree and number of variables in the algorithm. 
In order to compute the critical values, the following question should be solved. 
Let ${\mathscr A}$ be the set of homogeneous polynomials with $n$ variables of degree $d$  such that all nonzero coefficients are $1$, for example, $xyz+x^2y+z^3$. 
Then the symmetric group ${\frak S}_n$ acts on ${\mathscr A}$ by permuting the variables. 
The question is that how to get a set of representatives of the orbit space ${\mathscr A}/{\frak S}_n$.

In this paper, we also computed
\begin{itemize}
\item the list of $n$-variable monomials of degree $d$
\item the unique representation of a point in projective space ${\mathbb P}^n$ written by the coordinates of its affine cone ${\mathbb A}^{n+1}$ by fixing the first nonzero coordinate to be $1$
\item the degree and the list of coefficients of a multivariable polynomial $F$ which contains at least $2$ terms
\end{itemize}

Section \ref{MMaS} introduces the basic concepts of stability and moment maps, and how moment maps judge the stability in the complex setting. Then section \ref{MMoH} describes the general formulae of moment maps for hypersurfaces.

Theorem \ref{constrction} in section \ref{CtMM} gives the construction of the moment map via MAXIMA. Later in this section there is another construction using the general formulae (\ref{momentgeneral}) in Lemma \ref{lemmamomentgen}. With this construction we discuss the critical points of the square length function of the momoent map in section \ref{CPotSLoMM} and deduce that a real hypersurface is critical if and only if it is real critical, that is Theorem \ref{realcritical}. 

In section \ref{TDC} we compute critical points of the square length of the moment maps for plane cubics. In the last section \ref{N3D4} we compute some critical points for plane quartics.
 
\section{Moment Maps and Stability}\label{MMaS}

From now on, we assume that $G$ is a reductive affine algebraic group, and $X$ is a quasi-projective variety.
Denote the action of $G$ on $X$ by $g.x$ for $g\in G$ and $x \in X$. 
Let $\pi: L\rightarrow X$ be a line bundle on $X$. 
Suppose $L$ admits a $G$-action $g.l$ for $g\in G$ and $l\in L$. 
If $\pi(g.l)=g.\pi(l)$ for all $g\in G$ and  $l\in L$, we say this action is a $G$-{\bf linearization} of $L$ and $L$ is a {\bf linearized} $G$-bundle. 
For any line bundle $L$, there exists a number $n$ such that $L^{\otimes n}$ admits a $G$-linearization. 
\begin{definition} 
Let $L$ be a $G$-linearized line bundle on $X$ and $x\in X$.
\begin{itemize}
\item $x$ is called {\bf semi-stable} with respect to $L$ if there exists $m>0$ and $s\in \Gamma(X,L^m)^G$ (the $G$-invariant sections) such that $X_s=\{\, y\in X \, |\,  s(y)\neq 0\, \}$  is affine and contains $x$.
\item $x$ is called {\bf stable} with respect to $L$ if $x$ is semi-stable and additionally the isotropy group $G_x$ is finite and all orbits of $G$ in $X_s$ are closed where $X_s$ is defined in the definition above.
\item $x$ is called {\bf unstable} with respect to $L$ if it is not semi-stable.
\end{itemize}
\end{definition}
We shall denote the set of semi-stable (respectively stable, unstable) points by
$$X^{\rm ss}(L), \quad X^{\rm s}(L), \quad X^{\rm us}(L).$$
\begin{theorem}[Hilbert-Mumford Numerical Criteron of Stableness \cite{DolInv}] 
Let $\ {\rm Pic}^G(X)$ be the group of isomorphism classes of linearizations of the $G$-action. 
 For any ample
$L \in {\rm Pic}^G(X)$ on a complete $X$, there exists a real-valued function $M^L(x)$ such that 
$$X^{\rm ss} (L) = \{x \in X \ |\  M^L (x) \leq 0\},$$
$$X^s (L) = \{x \in X \ |\  M^L (x) < 0\}.$$
\end{theorem}

Consider the case when $X$ is projective with a specialized linearization $L$ such that the linear system $|L|$ determines an embedding $\iota_L: X\rightarrow {\mathbb P}^n$. Since $L$ is a linearization, one can show that the action of $G$ on $X$ is given by a representation $\rho :G\rightarrow {\rm GL}(n+1)$ while ${\rm GL}(n+1)$ acts linearly on ${\mathbb P}^n$. If we regard the analytic variety $X^{\rm an}$ as $X$, then $X$ is a symplectic manifold, whose symplectic form is endowed from ${\mathbb P}^n$, with the maximal compact subgroup (or the real form) $K$ of $G$ acts symplectically on $X$. We also assume that the image of $K$ under $\rho$ is contained in the unitary group ${\rm U}(n+1)$, the maximal compact subgroup of ${\rm GL}(n+1)$.

In general, let $X$ be a symplectic manofold with symplectic form $\omega$, let $K$ be a real Lie group acting on $X$ symplectically. Fix a point $x\in X$, define the map $K\rightarrow K.x$ by $g\mapsto g.x$ where $K.x$ is the orbit of $x$ in $X$.  This gives a set-theoretic exact sequence
$$K_x \longrightarrow K \longrightarrow K.x \subset X$$ where $K_x$ is the stabilizer of $x$. Because all arrows in this sequence are maps between differential manifolds, we have a diagram of tangent spaces at the identity $e\in K$ as 
$$T_e(K)\longrightarrow T_x(K.x)\subset T_x(X).$$ The tangent space $T_e(K)$ is the Lie algebra ${\frak g}:={\rm Lie}(K)$ of the group $K$. Denote the composition of these two maps by $\sharp:{\frak g}\rightarrow T_x(X)$. For $x\in X$, and any element $\xi\in {\frak g}$, we have a tangent vector ${\xi}^{\sharp}_x$. This defines a vector field ${\xi}^{\sharp}\in T(X)$.

\begin{definition}
A {\bf moment map} for the action of $K$ on $X$ is a smooth map $\Phi : X \rightarrow {\frak g}$ satisfying the following two properties:
\begin{itemize}
\item $\Phi$ is equivariant with respect to the action of $K$ on $X$ and the co-adjoint action of $K$
on ${\frak g}$ , 
\item for any $\xi \in {\frak g}$,we have $ \omega(\xi^{\sharp}, \cdot) = {\rm d}\, (\xi \circ \Phi)$.
\end{itemize}
\end{definition}
\begin{theorem}[See \cite{Ness}]\label{momentandML} Let $X$ be a projective variety with a reductive group $G$ action. Let $K$ be its maximal compact subgroup. Let $L$ be a linearization of the action and let $\Phi$ be the moment map constructed with respect to the representation of $K$ under the linear system $|L|$. Then	for any $x \in X$, the value $M^L (x)$ is equal to the signed distance from the origin to
the boundary of $\Phi(\, \overline{G.x}\, )$.
\end{theorem}

\section{Moment Maps of Hypersurfaces}\label{MMoH}

Let $V$ be complex vector space of dimension $m+1$. Let $G$ be a complex reductive algebraic group acting on $G$ by $G\times V\rightarrow V$, $(g,v)\mapsto g.v$. Let $K$ be the maximal compact subgroup of $G$. Assume that $V$ is equipped with a $K$-invariant hermitian form $\langle\ , \ \rangle$ and a $K$-invariant norm $\|v \|^2=\langle v,v\rangle$.  Let ${\mathbb P}(V)$ be the projective space of $V$, then $p_V:V\rightarrow {\mathbb P}(V)$, $(x_0,\cdots , x_m)\mapsto (x_0: \cdots : x_m)$ is the affine cone of ${\mathbb P}(V)$.

\begin{lemma}[See \cite{Ness}]\label{lemmamomentgen} For any $x\in {\mathbb P}(V)$, let $v \in V$ such that $p_V(v)=x$, define ${\frak m}(x)\in {\rm Lie}(K)^\vee$ such that for any $\alpha \in {\rm Lie}(K)$, 
\begin{equation}\label{momentgeneral}
{\frak m}(x)(\alpha)= \frac{1}{\|v\|^2} \left( \frac{{\rm d}}{{\rm d}t}\|{\rm e}^{t\alpha}.v\|^2 \right)_{t=0}
\end{equation}
then ${\frak m}:X\rightarrow {\rm Lie}(K)^{\vee}$ is a moment map with respect to $K$. 
\end{lemma}

Coming to hypersurfaces, we take
$$G={\rm SL}(n,{\mathbb C}),\quad K={\rm SU}(n,{\mathbb C}),\quad V={\rm Sym}^d({\mathbb C}^n)^\vee $$
where ${\rm Sym}^d({\mathbb C}^n)^\vee $ is the vector space of polynomials of homomgeneous degree $d$ with $n$ variables. Take $\{\, X_1,\ldots, X_n\, \}$ to be an orthogonoal basis of $({\mathbb C}^n)^\vee$. Denote $A=(a_{ij})_{n\times n}\in {\rm SL}(n, {\mathbb C})$ as the $n\times n$ matrix whose $(i,j)$-th entry is $a_{ij}$. Then the linear action of ${\rm SL}(n,{\mathbb C})$ on $V$ is defined by $A.X_i= \sum_{j=1}^n a_{ij}X_j.$ 
Let $$ {\mathscr W }=\{\, (\alpha_1,\ldots , \alpha_n)\in {\mathbb Z}^n \, | \, \alpha_1+\cdots +\alpha_n =d \, \}.$$
For any $\alpha=(\alpha_1,\ldots , \alpha_n) \in {\mathscr W}$, let $m^\alpha$ be the monomial $X_1^{\alpha_1}\cdots X_n^{\alpha_n}$. For any $f=\sum_{\alpha \in {\mathscr W}}c_\alpha m^{\alpha}$ where $c_\alpha \in {\mathbb C}$, then the action of ${\rm SL}(n,{\mathbb C})$ on ${\rm Sym}^d({\mathbb C}^n)^\vee$ is just the substitution of the action of action on $({\mathbb C}^n)^\vee$. That is, for any $A\in {\rm SL}(n, {\mathbb C})$, we have 
$$A.f=\sum_{(\alpha_1,\ldots , \alpha_n)\in {\mathscr W}} \left( c_{\alpha}\cdot \prod_{i=1}^n \left(\sum_{j=1}^n a_{ij}X_j\right)^{\alpha_i}\right).$$

According to \cite{Ness}, an ${\rm SU}(n)$-invariant hermitian product is defined by 
\begin{itemize}
\item The monomials $\{ m^{\alpha}\}_{\alpha \in {\mathscr W}}$ form an orthogonal basis of ${\rm Sym}^d({\mathbb C}^n)^\vee$.
\item $\displaystyle{ \|m^{\alpha}\|^2 = \frac{\alpha_1！\cdots \alpha_n!}{d!}}$.
\end{itemize} Denote this hermitian product by $\langle \, , \, \rangle$.

For any $f=\sum_{\alpha\in {\mathscr W}}c_{\alpha}m^{\alpha}\in {\rm Sym}^d({\mathbb C}^n)^\vee$, let $H(f)$ be the $n\times n$ Hermitian matrix whose $(i,j)$-th entry is
\begin{equation}H(f)^j_i = \frac{1}{\|f\|^2}\cdot \frac{1}{d}\, \left\langle \frac{\partial f}{\partial x_j}\, , \frac{\partial f}{\partial x_i} \right\rangle. \end{equation}
\begin{lemma} Let $I$ be the $n\times n$ identity matrix. The map
\begin{equation} \begin{array}{cccc}
{\frak m}^*: & {\mathbb P}({\rm Sym}^d({\mathbb C}^n)^\vee)& \longrightarrow & \sqrt{-1}\ {\frak{su}}(n)\\
& f & \longmapsto & \displaystyle{2\left(H(f)-\frac{d}{n}I\right)}
\end{array}\end{equation}
is the dual of the moment map. Since ${\frak s \frak u}(n)$ is semisimple, the Killing form identifies ${\frak{su}}(n)$ with its dual. We confuse ${\frak m}$ with ${\frak m}^*$ here.
\end{lemma}

\section{Constructing the Moment Map}\label{CtMM}
We will take $n=3$ and $d=3$ as an example, but the code we write here can compute any other case by simply changing $n$ and $d$ in the claiming part at the very beginning of the code.

Not all objects and algorithms we need are prepared in MAXIMA, but some of them are very elementary and important. So it is not harmful to display them in details. 

We start with constructing the list of monomials for a specific degree $d$ and number of variables $n$. Though MAXIMA has such a function in the package {\tt affine}, the outputs make trouble while one takes substitution. The list we make here is univeral. The idea is to expand $(X_1+\cdots +X_n)^d$ and take the terms. First we claim $n,d$ and the name of variables.
\begin{flushleft}

\begin{verbatim}
block(n:3, d:3, X:[x,y,z],Y:create_list(i=1,i,X), display(n,d,X));
\end{verbatim}

\end{flushleft}
Then define the function of monomials by
\begin{flushleft}

\begin{verbatim}
M(k):=sort(create_list(part( 
     expand(lsum(i,i,X)^k),j)/subst(Y,part(expand(lsum(i,i,X)^k),j)),
     j,1,binomial(n+k-1,k)));
\end{verbatim}

\end{flushleft}
Note that this function can also be defined as a function of both $k$ and $X$. Here we do not change the variables, so we let $M$ only depend on the degree $k$. One sees that $\tt M(2)$ is $[{{x}^{2}},x y,{{y}^{2}},x z,y z,{{z}^{2}}]$ and $\tt M(d)$ is \[[{{x}^{3}},{{x}^{2}} y,x {{y}^{2}},{{y}^{3}},{{x}^{2}} z,x y z,{{y}^{2}} z,x {{z}^{2}},y {{z}^{2}},{{z}^{3}}].\]

Furthurmore, we define the weight and degree of a monomial as

\begin{flushleft}

\begin{verbatim}
Wt(F):=create_list(hipow(F,i),i,X);
\end{verbatim}

\begin{verbatim}
deg(F):=lsum(i,i,Wt(F));
\end{verbatim}

and if $F$ is a polynomial with at least $2$ terms, then the degree of $F$ is
\begin{verbatim}
Deg(F):= lmax(create_list(deg(part(F,i)),i,1,length(F)));
\end{verbatim}

\end{flushleft}

With these preparation, we construct the moment map of a general complex polynomial
\begin{flushleft}

\begin{verbatim}
C:lsum((Re[Wt(j)[1],Wt(j)[2],Wt(j)[3]]
        +%i*Im[Wt(j)[1],Wt(j)[2],Wt(j)[3]])*j,j,M);
\end{verbatim}

$\left( {\rm i} {{\mathit{Im}}_{0,0,3}}+{{\mathit{Re}}_{0,0,3}}\right)  {{z}^{3}}+\left( {\rm i} {{\mathit{Im}}_{0,1,2}}+{{\mathit{Re}}_{0,1,2}}\right)  y {{z}^{2}}+\left( {\rm i} {{\mathit{Im}}_{1,0,2}}+{{\mathit{Re}}_{1,0,2}}\right)  x {{z}^{2}}+\left( {\rm i} {{\mathit{Im}}_{0,2,1}}+{{\mathit{Re}}_{0,2,1}}\right)  {{y}^{2}} z+\left( {\rm i} {{\mathit{Im}}_{1,1,1}}+{{\mathit{Re}}_{1,1,1}}\right)  x y z+\left( {\rm i} {{\mathit{Im}}_{2,0,1}}+{{\mathit{Re}}_{2,0,1}}\right)  {{x}^{2}} z+\left( {\rm i} {{\mathit{Im}}_{0,3,0}}+{{\mathit{Re}}_{0,3,0}}\right)  {{y}^{3}}+\left( {\rm i} {{\mathit{Im}}_{1,2,0}}+{{\mathit{Re}}_{1,2,0}}\right)  x {{y}^{2}}+\left( {\rm i} {{\mathit{Im}}_{2,1,0}}+{{\mathit{Re}}_{2,1,0}}\right)  {{x}^{2}} y+\left( {\rm i} {{\mathit{Im}}_{3,0,0}}+{{\mathit{Re}}_{3,0,0}}\right)  {{x}^{3}}$

\end{flushleft}

Now we define the ${\rm SU}(n)$-invariant hermitian product on ${\rm Sym}^d({\mathbb C}^n)^\vee$. For a general homoegeneous polynomial $F$, define the weighted coefficient list by

\begin{flushleft}

\begin{verbatim}
wt(F):=product(Wt(F)[i]!,i,1,length(Wt(F)))/deg(F)!;

WtCoe(F):=subst(create_list(k=1,k,X),create_list(
   sum(sqrt(wt(part(F,j)))*part(F,j)*kron_delta(Wt(part(F,j)),Wt(i)),
   j,1,length(F)),i,M(Deg(F))));
\end{verbatim}

\end{flushleft}
Then the hermitian product $\langle F,G \rangle$ is
\begin{flushleft}
\begin{verbatim}
conjugate(WtCoe(F)).WtCoe(G);
\end{verbatim}
\end{flushleft}

\begin{theorem}\label{constrction} Let $C$ be the general complex polynomial, construct
\begin{flushleft}
\begin{verbatim}
block(
    LCC:conjugate(WtCoe(C)).WtCoe(C),
    G:create_list(diff(C,i),i,X),
    PP[i,j]:=conjugate(WtCoe(G[i])).WtCoe(G[j]),
    HH:expand(genmatrix(PP,n,n)/d),
    m:ratsimp(expand(2*(HH/LCC-diagmatrix(n,d/n)))),
    Lm:ratsimp(realpart(mat_trace(m.m)))
);
\end{verbatim}
\end{flushleft}
Then $m$ is the moment matrix of $C$, the expression $Lm$ is the square length $\|{ m}\|^2$ of $ m$. Let $R$ be the general real polynomial

\begin{flushleft}

\begin{verbatim}
R:lsum(a[Wt(j)[1],Wt(j)[2],Wt(j)[3]]*j,j,M(d));
\end{verbatim}
${{a}_{0,0,3}} {{z}^{3}}+{{a}_{0,1,2}} y {{z}^{2}}+{{a}_{1,0,2}} x {{z}^{2}}+{{a}_{0,2,1}} {{y}^{2}} z+{{a}_{1,1,1}} x y z+{{a}_{2,0,1}} {{x}^{2}} z+{{a}_{0,3,0}} {{y}^{3}}+{{a}_{1,2,0}} x {{y}^{2}}+{{a}_{2,1,0}} {{x}^{2}} y+{{a}_{3,0,0}} {{x}^{3}}$
\end{flushleft}
then by replacing $C$ into $R$, the construction above gives the moment matrix $m$ of $R$ and the square length $\|{m}\|^2$ of $R$.
\end{theorem}

For easiness, we only consider real cubic polynomials from now on. Let $A$ be a random real cubic homogeneous polynomial in ${\rm Sym}^d({\mathbb R}^n)^\vee$ which has at least 2 terms (our weight function only defines for such a polynomial). Define the coefficient list by
\begin{flushleft}
\begin{verbatim}
Coe(F):=subst(create_list(k=1,k,X),create_list(
     sum(part(F,j)*kron_delta(Wt(part(F,j)),Wt(i)),
     j,1,length(F)),i,M(Deg(F))));
\end{verbatim}
\end{flushleft}
Then the list

\begin{flushleft}
\begin{verbatim}
SB(A):=create_list(Coe(R)[i]=Coe(A)[i],i,1,length(Coe(R)));
\end{verbatim}
\end{flushleft}
substitudes the coefficients of $A$ into $R$. So the function

\begin{flushleft}
\begin{verbatim}
Moment(A):=subst(SB(A),m);
\end{verbatim}
\end{flushleft} 
gives the moment matrix of $A$.

\begin{example} The moment matrix of $x^3+y^3$ is
\[\begin{pmatrix}1 & 0 & 0\cr 0 & 1 & 0\cr 0 & 0 & -2\end{pmatrix}\]
\end{example}

There is another method to construct the moment matrix of a given polynomial using (\ref{momentgeneral}). This method is much easier to be generalized. For a matrix $A\in {\frak gl}(n+1)$ and a polynomial $f(x,y,z)$ where $X$ is the list $[x,y,z]$, we define the action ${\rm e}^{tA}.f$ by
\begin{flushleft}
\begin{verbatim}
Act(A,f):=expand(subst(
  create_list(X[i]=(ratsimp (mat_function (exp, t*A)).X)[i][1],
  i,1,d),f));
\end{verbatim}
\end{flushleft}
Then the $(i,j)$-th entry of the moment matrix of $R$ is
\begin{flushleft}
\begin{verbatim}
for i=1 thru d do for j=1 thru d do
ratsimp(subst([t=0],diff(SqL(Act(ematrix(d,d,1,i,j),R))/2,t)));
\end{verbatim}
\end{flushleft}

\newpage \section{Stability of Plane Cubics} \label{SOPC}
\subsection{Critical Points of the Square Length of the Moment Map}\label{CPotSLoMM}

Given a moment map ${\frak m}: V\rightarrow {\rm Lie}(K)^\vee$ where $V$ is a vector space with a $G$-action, and $K$ is the maximal compact subgroup of $G$, it is very important to study the critical points of $\|{\frak m}\|^2$. 
\begin{theorem}\label{criticalorbit} Let $\|{\frak m}\|^2 \big|_{G.x}$ denote the restriction of $\|{\frak m}\|^2$ on the $G$-orbit of $x$ in ${\mathbb P}(V)$. If ${\rm d}\|{\frak m}\|^2(x)=0$, then
\begin{enumerate}
\item $\|{\frak m}\|^2 \big|_{G.x}$ attains its minimum value at $x$,
\item $\|{\frak m}\|^2 \big|_{G.x}$ attains its minimum on a unique $K$-orbit.
\end{enumerate}
\end{theorem}

For hypersurfaces, we have $$\|{\frak m}\|^2= {\rm Re}\, (\, {\rm Tr}\, (\, {\frak m}.{\frak m}\, )\, ).$$
As mentioned before, we only compute the real critical values. For a general complex polynomial $C$, the gradient of $C$ is the vector $$ \Theta =\left(\frac{\partial\| {\frak m}(C)\|^2}{Re_{i,j,k}}\, , \, \frac{\partial\| {\frak m}(C)\|^2}{Im_{i,j,k}}\right)_{(i,j,k)\in {\mathscr W}}.$$ 

\begin{theorem}\label{realcritical}
If $R$ is a polynomial defined over the real numbers ${\mathbb R}$, then $R$ is a critical value of $\|{\frak m}\|^2$ if and only if it is a real critical value of $\|{\frak m}\|^2$.
\end{theorem}
{\bf Proof}\quad First we have the complex polynomial

\begin{verbatim}
C:lsum((a[Wt(j)[1],Wt(j)[2],Wt(j)[3]]+
        %i*b[Wt(j)[1],Wt(j)[2],Wt(j)[3]])*j,j,M(d));
\end{verbatim}
\begin{align*}&
\left( {{a}_{0,0,3}}+{\rm i} {{b}_{0,0,3}}\right)  {{z}^{3}}+\left( {{a}_{0,1,2}}+{\rm i}\,  {{b}_{0,1,2}}\right)  y {{z}^{2}}+\left( {{a}_{1,0,2}}+{\rm i}\,  {{b}_{1,0,2}}\right)  x {{z}^{2}}+ \\ & \left( {{a}_{0,2,1}}+{\rm i}\,  {{b}_{0,2,1}}\right)  {{y}^{2}} z+\left( {{a}_{1,1,1}}+{\rm i}\,  {{b}_{1,1,1}}\right)  x y z+\left( {{a}_{2,0,1}}+{\rm i}\,  {{b}_{2,0,1}}\right)  {{x}^{2}} z+ \\ & \left( {{a}_{0,3,0}}+{\rm i}\,  {{b}_{0,3,0}}\right)  {{y}^{3}}+\left( {{a}_{1,2,0}}+{\rm i}\,  {{b}_{1,2,0}}\right)  x {{y}^{2}}+\left( {{a}_{2,1,0}}+{\rm i}\,  {{b}_{2,1,0}}\right)  {{x}^{2}} y+ \\ & \left( {{a}_{3,0,0}}+{\rm i}\,  {{b}_{3,0,0}}\right)  {{x}^{3}}.\end{align*}

Since both the functions ``{\tt realpart()}" and ``{\tt imagpart()}" can be applied to lists as
\begin{verbatim}
realpart(Coe(C));
\end{verbatim}
$$
[\, {{a}_{3,0,0}},{{a}_{2,1,0}},{{a}_{1,2,0}},{{a}_{0,3,0}},{{a}_{2,0,1}},{{a}_{1,1,1}},{{a}_{0,2,1}},{{a}_{1,0,2}},{{a}_{0,1,2}},{{a}_{0,0,3}}\, ]$$

\begin{verbatim}
imagpart(Coe(C));
\end{verbatim}
$$
[\, {{b}_{3,0,0}},{{b}_{2,1,0}},{{b}_{1,2,0}},{{b}_{0,3,0}},{{b}_{2,0,1}},{{b}_{1,1,1}},{{b}_{0,2,1}},{{b}_{1,0,2}},{{b}_{0,1,2}},{{b}_{0,0,3}}\, ]$$
we can define the list
\begin{verbatim}
CCoe(C):=append(realpart(Coe(C)),imagpart(Coe(C)));
\end{verbatim}
such that 
\begin{verbatim}
CCoe(C);
\end{verbatim}
\begin{align*}
&[\, {{a}_{3,0,0}},{{a}_{2,1,0}},{{a}_{1,2,0}},{{a}_{0,3,0}},{{a}_{2,0,1}},{{a}_{1,1,1}},{{a}_{0,2,1}},{{a}_{1,0,2}},{{a}_{0,1,2}},{{a}_{0,0,3}}, \\ &\quad {{b}_{3,0,0}},{{b}_{2,1,0}},{{b}_{1,2,0}},{{b}_{0,3,0}},{{b}_{2,0,1}},{{b}_{1,1,1}},{{b}_{0,2,1}},{{b}_{1,0,2}},{{b}_{0,1,2}},{{b}_{0,0,3}}\, ].\end{align*}

Then we have the complex version of the moment matrix
\begin{verbatim}
block(
    cLCC:SqL(C),
    cG:create_list(diff(C,i),i,X),
    cPP[i,j]:=conjugate(WtCoe(cG[i])).WtCoe(cG[j]),
    cHH:expand(genmatrix(cPP,n,n)/d),
    cm:ratsimp(expand(2*(cHH/cLCC-diagmatrix(n,d/n)))),
    cLm:ratsimp(realpart(mat_trace(cm.cm))),
    cGRD: create_list(diff(cLm,i),i,CCoe(C)),
    0
);
\end{verbatim}
For the list
\begin{verbatim}
CCoe(R);
\end{verbatim}
\begin{align*}
&[{{a}_{3,0,0}},{{a}_{2,1,0}},{{a}_{1,2,0}},{{a}_{0,3,0}},{{a}_{2,0,1}},{{a}_{1,1,1}},{{a}_{0,2,1}},{{a}_{1,0,2}},{{a}_{0,1,2}},{{a}_{0,0,3}},\\ &\quad  0,0,0,0,0,0,0,0,0,0]\end{align*}
we define a more general substitution function 
\begin{verbatim}
Sb(A,B):=create_list(A[i]=B[i],i,1,length(A));
\end{verbatim}
such that 
\begin{verbatim}
Sb(CCoe(C),CCoe(R));
\end{verbatim}
\begin{align*}
& [{{a}_{3,0,0}}={{a}_{3,0,0}},{{a}_{2,1,0}}={{a}_{2,1,0}},{{a}_{1,2,0}}={{a}_{1,2,0}},{{a}_{0,3,0}}={{a}_{0,3,0}},{{a}_{2,0,1}}={{a}_{2,0,1}}, \\ & {{a}_{1,1,1}}={{a}_{1,1,1}},{{a}_{0,2,1}}={{a}_{0,2,1}},{{a}_{1,0,2}}={{a}_{1,0,2}},{{a}_{0,1,2}}={{a}_{0,1,2}},{{a}_{0,0,3}}={{a}_{0,0,3}},\\ & {{b}_{3,0,0}}=0,{{b}_{2,1,0}}=0,{{b}_{1,2,0}}=0,{{b}_{0,3,0}}=0,{{b}_{2,0,1}}=0,{{b}_{1,1,1}}=0,{{b}_{0,2,1}}=0,\\ & {{b}_{1,0,2}}=0,{{b}_{0,1,2}}=0,{{b}_{0,0,3}}=0]\end{align*}
The theorem is proved by
\begin{verbatim}
subst(Sb(CCoe(C),CCoe(R)),
      create_list(diff(cLm,i),i,imagpart(Coe(C))));
\end{verbatim}
$$
[\, 0,0,0,0,0,0,0,0,0,0\, ]$$

\begin{flushleft}
$\blacksquare$
\end{flushleft}

For a general real polynomial $A\in {\rm Sym}^d({\mathbb C}^n)^\vee$, define the real gradient of $A$ by
\begin{flushleft}

\begin{verbatim}
block(
    Lm:ratsimp(realpart(mat_trace(m.m))),
    GRD: create_list(diff(Lm,i),i,Coe(R))
);

Grad(A):=subst(SB(A),GRD);
\end{verbatim}
\end{flushleft}

Next we wish to find solutions for the equation ${\rm Grad}(A)=0$. The computation is too heavy for a general laptop. On the other hand, \cite{Ness} gives a list of critical points of plane cubics, it seems that they all have much less terms than the general polynomial $R$. Also, the algorithm can be simplified.

\begin{lemma} A polynomial $f\in {\rm Sym}^d({\mathbb C}^n)^\vee$ is a critical point of $\|{\frak m}\|^2$ if and only if ${\frak m}^*(f)$ is fixed by a diagonal on parameter subgroup $e^{t\alpha}$ for some diagonal $\alpha \in {\frak{su}}(n)$ and ${\frak m}^*(f)$ is a multiple of ${\alpha}$. Thus, if $f$ is a critical point, then ${\frak m}^*(f)$ is diagonal.
\end{lemma} 

In the next section we first find polynomials whose moment matrices are diagonal. 

\subsection{Representatives of Symmetric Orbits, Diagonal Moment Matrices and Critical Points of Plane Cubics}\label{TDC}

It is easy to compute the moment matrices of the monomials.
\begin{theorem}[see \cite{Ness}] If $f\in {\rm Sym}^d({\mathbb C}^n)^\vee$ is a monomial, then ${\frak m}^*(f)$ is diagonal and $f$ is a critical point of $\|{\frak m}\|^2$.
\end{theorem}

Regard $f$ as an element in $f\in{\mathbb{P}}( {\rm Sym}^d({\mathbb C}^n)^\vee)$, we can assume that all coefficients of $f$ are equal to $1$. For example, if we have computed that the moment matrix of $x^3$ is diagonal, then it is easy to imagine that so is $y^3$ and $z^3$. Because the symmetric group ${\frak S}_n$ acts on ${\mathbb C}^n$ by permuting the coordinates, it also acts on ${\mathbb{P}}(f\in {\rm Sym}^d({\mathbb C}^n)^\vee)$ by the substitution of the permutation. Thus one can imagine that if $f$ and $g$ are in the same ${\frak S}_n$-orbit, then one of them has a diagonal moment matrix implies that the other also has a diagonal moment matrix. Suppose we compute the moment matrix of the polynomial $ax^3+bx^2y$, then we do not have to compute $cy^3+dy^2z$ because $x^3+x^2y$ and $y^3+y^2z$ are in the same ${\frak S}_n$-orbit.

Let ${\mathscr T}^{n,d}$($={\mathscr T}$ for short) denote the set of representatives of ${\frak S}_n$-orbits of  for which each element is the minimum in its orbit according to the lexicograpic order, and let ${\mathscr T}_m\subset {\mathscr T}$ be the subset of elements which have exactly $m$-terms.

\begin{example} Let $n=d=3$. Assume that $x\prec y \prec z$, then ${\mathscr T}_1=\{x^3,x^2y,xyz\}$.
\end{example}

The function ``{\tt orbit}" in Maxima computes the ${\frak S}_n$-orbit of a polynomial $f$ with respect to a list $X$ of length $n$ (in fact ${\frak S}_n$ acts on the set of $X$). In order to compute ${\mathscr T}_m$, we construct the equivalent classes with respect to the relation $f\sim g$ if and only if $f$ lies in the ${\frak S}_n$-orbit of $g$. The stratification we make for ${\mathscr T}$ is for the purpose to deduce the computing amount of a laptop for a single cell.

First, we need the power set of the set of monomials.

\begin{flushleft}
\begin{verbatim}
block(
     B:binomial(n+d-1,d),
     N: setify(M(d)),
     PN:listify(powerset(N)));
\end{verbatim}
\end{flushleft}
\begin{proposition}
Define the ${\frak S}_n$-conjugate equivalent relation by 

\begin{flushleft}
\begin{verbatim}
ff(g,h):= member(g,orbit(h,[x,y,z]));
\end{verbatim}
\end{flushleft}

Then the function
\begin{flushleft}
\begin{verbatim}
U(m):=create_list(t[1],t,sort(full_listify(
      equiv_classes(setify(create_list(
      lsum(i,i,listify(t)),
      t,listify(subset(setify(PN),
                lambda([x],is(cardinality(x)=m)
       ))))),ff))));
\end{verbatim}
\end{flushleft}
outputs the list of the set ${\mathscr T}_m$.
\end{proposition}
\begin{example} For $n=d=3$, we have that ${\mathscr T}_2$ is

\begin{flushleft}
$\{\,  {{x}^{2}} y+{{x}^{3}},x {{y}^{2}}+{{x}^{3}},x {{y}^{2}}+{{x}^{2}} y,{{y}^{3}}+{{x}^{3}},{{x}^{2}} z+{{x}^{2}} y,{{x}^{2}} z+x {{y}^{2}},{{x}^{2}} z+{{y}^{3}},x y z+{{x}^{3}},x y z+{{x}^{2}} y,{{y}^{2}} z+{{x}^{2}} z\}$
\end{flushleft} and ${\mathscr T}_3$ is
\begin{flushleft}
$\{\,  x {{y}^{2}}+{{x}^{2}} y+{{x}^{3}},{{y}^{3}}+{{x}^{2}} y+{{x}^{3}},{{x}^{2}} z+{{x}^{2}} y+{{x}^{3}},{{x}^{2}} z+x {{y}^{2}}+{{x}^{3}},{{x}^{2}} z+x {{y}^{2}}+{{x}^{2}} y,{{x}^{2}} z+{{y}^{3}}+{{x}^{3}},{{x}^{2}} z+{{y}^{3}}+{{x}^{2}} y,{{x}^{2}} z+{{y}^{3}}+x {{y}^{2}},x y z+{{x}^{2}} y+{{x}^{3}},x y z+x {{y}^{2}}+{{x}^{3}},x y z+x {{y}^{2}}+{{x}^{2}} y,x y z+{{y}^{3}}+{{x}^{3}},x y z+{{x}^{2}} z+{{x}^{2}} y,x y z+{{x}^{2}} z+x {{y}^{2}},x y z+{{x}^{2}} z+{{y}^{3}},{{y}^{2}} z+{{x}^{2}} z+{{x}^{3}},{{y}^{2}} z+{{x}^{2}} z+{{x}^{2}} y,{{y}^{2}} z+x y z+{{x}^{2}} z,x {{z}^{2}}+x {{y}^{2}}+{{x}^{3}},x {{z}^{2}}+x {{y}^{2}}+{{x}^{2}} y,x {{z}^{2}}+{{y}^{3}}+{{x}^{3}},x {{z}^{2}}+{{y}^{3}}+{{x}^{2}} y,x {{z}^{2}}+{{x}^{2}} z+{{y}^{3}},x {{z}^{2}}+{{y}^{2}} z+{{x}^{2}} y,{{z}^{3}}+{{y}^{3}}+{{x}^{3}}\, \}.$
\end{flushleft}
\end{example}

According to \cite{Ness}, the hypersurfaces containing less than $n$ variables will give critical points coming from a lower dimensional case. So we only need to compute moment matrix of a hypersurface such that all $x,y,z$ are in its expression. The function `{\tt freeof}' outputs true if an expression is not in another one. So the function

\begin{flushleft}
\begin{verbatim}
NF(x,A):= not(freeof(x,A));
LNF(x,A):= freeof(false,create_list(NF(i,A),i,x));
fff(h):=LNF(X,h);
V(n):=sublist(U(n),fff);
\end{verbatim}
\end{flushleft}
will select elements in ${\mathscr T}_m$ such that all $x,y,z$ are in the expressions. 

We will compute the moment matrices of 
\begin{flushleft}
\begin{verbatim}
T(n):=sort(create_list(sum(
      b[i]*part(j,i),i,1,n-1)+part(j,n),j,V(n)));
\end{verbatim}
\end{flushleft}
and pick up those who have diagonal moment matrices.

For $f\in{\mathbb{P}}( {\rm Sym}^d({\mathbb C}^n)^\vee)$, let $l(f)$ be the number of terms of $f$.
\begin{theorem}\label{diag} For $n=d=3$, if $l(f)\geq 5$, then the moment matrix ${\frak m}(f)$ is not diagonal. For $2\leq l(f) \leq 4$, the polynomials of the form 
\begin{align*}
&{{b}_{1}} {{x}^{2}} z+x {{y}^{2}},\ {{b}_{1}} {{x}^{2}} z+{{y}^{3}},\ {{b}_{1}} x y z+{{x}^{3}},\ {{b}_{1}} {{y}^{2}} z+{{x}^{2}} z, \\&{{b}_{1}} x y z+{{b}_{2}} {{y}^{3}}+{{x}^{3}},\ {{b}_{1}} x {{z}^{2}}+{{b}_{2}} x {{y}^{2}}+{{x}^{3}}, \\ &{{b}_{1}} x {{z}^{2}}+{{b}_{2}} {{y}^{3}}+{{x}^{3}},\ {{b}_{1}} x {{z}^{2}}+{{b}_{2}} {{y}^{3}}+{{x}^{2}} y,\\ &{{b}_{1}} x {{z}^{2}}+{{b}_{2}} {{y}^{2}} z+{{x}^{2}} y,\ {{b}_{1}} {{z}^{3}}+{{b}_{2}} {{y}^{3}}+{{x}^{3}}, \\ &{{b}_{1}} {{z}^{3}}+{{b}_{2}} x y z+{{b}_{3}} {{y}^{3}}+{{x}^{3}}
\end{align*}
gives all the polynomials who have diagonal matrix up to a scalar multiplication.
\end{theorem}

\paragraph{Remark of Theorem  \ref{diag}} In order to prove Theorem  \ref{diag}, one computes the moment matrices at one hand, but also have to check that for those whose moment matrices are not diagonal as a general expression, the non-diagonal entries are all zero only if one of their coefficients vanishes. 

\begin{example}\label{criticalcubic} Using ``{\tt algsys}", we solve the critical points for $l(f)\leq 3$. obvious isomorphism (a scalar product on one of the coordinates), they are
\begin{align*}
& {{y}^{2}} z+{{x}^{2}} z,\ {{x}^{2}} z+x {{y}^{2}}, \\
& {{z}^{3}}+{{y}^{3}}+{{x}^{3}}, \ x {{z}^{2}}+{{y}^{2}} z+{{x}^{2}} y ,\\ & \sqrt{2} x {{z}^{2}}+\frac{{{y}^{3}}}{\sqrt{3}}+{{x}^{2}} y ,\ 3 x {{z}^{2}}+\sqrt{2} {{y}^{3}}+{{x}^{3}}.
\end{align*}

We will introduce the precise computation in the next section. Here
we use Mathematica to draw the affine cones of these curves over the real field ${\mathbb R}$ (see Figure \ref{affineconecubic}).

\end{example}

According to \cite{Ness}, the family $a({{z}^{3}}+{{y}^{3}}+{{x}^{3}})+bxyz$ are critical points, and ${{z}^{3}}+{{y}^{3}}+{{x}^{3}}$ is contained in this family. In order to get this whole family of critical points, we consider hypersurfaces of type ${{b}_{1}} {{z}^{3}}+{{b}_{2}} x y z+{{b}_{3}} {{y}^{3}}+{{x}^{3}}$. For easiness, we consider ${{b}_{3}} {{z}^{3}}+ x y z+{{b}_{2}} {{y}^{3}}+{{b}_{1}}{{x}^{3}}$ which defines the same type of hypersurfaces since all coefficients are not zero. 
\begin{lemma}[See \cite{Ness}]\label{fixed} Let $f\in{\mathbb{P}}( {\rm Sym}^d({\mathbb C}^n)^\vee)$ defines a hypersurface, then $f$ is a critical point of $\|{\frak m}\|^2$ if and only if ${\frak m}(f)$ is diagonal and ${\rm e}^{{\frak m}(f)}$ fixes $f$ in ${\mathbb{P}}( {\rm Sym}^d({\mathbb C}^n)^\vee)$.
\end{lemma}

\begin{proposition}
Let $F$ be the real polynomial ${{b}_{3}} {{z}^{3}}+ x y z+{{b}_{2}} {{y}^{3}}+{{b}_{1}}{{x}^{3}}$ such that $b_1b_2b_3\neq 0$. Then $F$ is a critical value of $\|{\frak m}\|^2$ if and only if $b_1^2=b_2^2=b_3^2$. That is, $F$ is obviously isomorphic to the family $\lambda xyz + \mu (x^3+y^3+z^3)$.
\end{proposition} 
{\bf Proof}\quad  By lemma \ref{fixed}, we compute ${\rm e}^{{\frak m}(f)}.f$, and output the coefficient vector in ${\mathbb{P}}( {\rm Sym}^d({\mathbb C}^n)^\vee)$.

\begin{flushleft}
\begin{verbatim}
subst([t=1],Pj(Coe(Act(Moment(F),F))));
\end{verbatim}
where $\tt Pj()$ is defined as in \ref{projectivespace}, we have
\begin{align*}
 & \left[1,0,0, \frac{{{e}^{\frac{108 {{b}_{2}^{2}}}{6 {{b}_{3}^{2}}+6 {{b}_{2}^{2}}+6 {{b}_{1}^{2}}+1}-\frac{108 {{b}_{1}^{2}}}{6 {{b}_{3}^{2}}+6 {{b}_{2}^{2}}+6 {{b}_{1}^{2}}+1}}} {{b}_{2}}}{{{b}_{1}}},0, \right. \\& \ \\ & \frac{{{e}^{-\frac{72 {{b}_{1}^{2}}}{6 {{b}_{3}^{2}}+6 {{b}_{2}^{2}}+6 {{b}_{1}^{2}}+1}+\frac{36 {{b}_{2}^{2}}}{6 {{b}_{3}^{2}}+6 {{b}_{2}^{2}}+6 {{b}_{1}^{2}}+1}+\frac{36 {{b}_{3}^{2}}}{6 {{b}_{3}^{2}}+6 {{b}_{2}^{2}}+6 {{b}_{1}^{2}}+1}}}}{{{b}_{1}}},0,0,0, 
\\& \ \\ &  \left. \frac{{{e}^{\frac{108 {{b}_{3}^{2}}}{6 {{b}_{3}^{2}}+6 {{b}_{2}^{2}}+6 {{b}_{1}^{2}}+1}-\frac{108 {{b}_{1}^{2}}}{6 {{b}_{3}^{2}}+6 {{b}_{2}^{2}}+6 {{b}_{1}^{2}}+1}}} {{b}_{3}}}{{{b}_{1}}}\ \right]
\end{align*}
\end{flushleft}

Comparing this vector with the coefficient vector of $F$ which is 
$$\left[1,0,0,\frac{{{b}_{2}}}{{{b}_{1}}},0,\frac{1}{{{b}_{1}}},0,0,0,\frac{{{b}_{3}}}{{{b}_{1}}}\right],$$
one sees that $b_1^2=b_2^2=b_3^2$.
\newpage \section{Some Unstable Plane Quartics} \label{SUPQ}

Now let us come to the case when $n=3$ and $d=4$. We only compute unstable curves with 2 or 3 terms. 
The moment matrix of 
\begin{align*}
R=& {{a}_{0,0,4}} {{z}^{4}}+{{a}_{0,1,3}} y {{z}^{3}}+{{a}_{1,0,3}} x {{z}^{3}}+{{a}_{0,2,2}} {{y}^{2}} {{z}^{2}}+{{a}_{1,1,2}} x y {{z}^{2}}+{{a}_{2,0,2}} {{x}^{2}} {{z}^{2}}+ \\ & {{a}_{0,3,1}} {{y}^{3}} z+{{a}_{1,2,1}} x {{y}^{2}} z+{{a}_{2,1,1}} {{x}^{2}} y z+{{a}_{3,0,1}} {{x}^{3}} z+{{a}_{0,4,0}} {{y}^{4}}+{{a}_{1,3,0}} x {{y}^{3}}+ \\ & {{a}_{2,2,0}} {{x}^{2}} {{y}^{2}}+{{a}_{3,1,0}} {{x}^{3}} y+{{a}_{4,0,0}} {{x}^{4}}
\end{align*} is the symmetric matrix
$$\frac{1}{r}\cdot 
\begin{pmatrix}{{r}_{1,1}} & {{r}_{1,2}} & {{r}_{1,3}}\cr {{r}_{2,1}} & {{r}_{2,2}} & {{r}_{2,3}}\cr {{r}_{3,1}} & {{r}_{3,2}} & {{r}_{3,3}}\end{pmatrix}$$

where

\begin{align*}
 r=& 36 {{a}_{4,0,0}^{2}}+9 {{a}_{3,1,0}^{2}}+9 {{a}_{3,0,1}^{2}}+6 {{a}_{2,2,0}^{2}}+3 {{a}_{2,1,1}^{2}}+6 {{a}_{2,0,2}^{2}}+  9 {{a}_{1,3,0}^{2}}+ \\& 3 {{a}_{1,2,1}^{2}}+3 {{a}_{1,1,2}^{2}}+9 {{a}_{1,0,3}^{2}}+36 {{a}_{0,4,0}^{2}}+9 {{a}_{0,3,1}^{2}}+6 {{a}_{0,2,2}^{2}}+9 {{a}_{0,1,3}^{2}}+36 {{a}_{0,0,4}^{2}}\\
 {{r}_{1,1}}=&-96 {{a}_{0,0,4}^{2}}-24 {{a}_{0,1,3}^{2}}-16 {{a}_{0,2,2}^{2}}-24 {{a}_{0,3,1}^{2}}-96 {{a}_{0,4,0}^{2}}-6 {{a}_{1,0,3}^{2}}\\& -2 {{a}_{1,1,2}^{2}}- 2 {{a}_{1,2,1}^{2}}-6 {{a}_{1,3,0}^{2}}+8 {{a}_{2,0,2}^{2}}+4 {{a}_{2,1,1}^{2}}+8 {{a}_{2,2,0}^{2}}+30 {{a}_{3,0,1}^{2}}+\\& 30 {{a}_{3,1,0}^{2}}+192 {{a}_{4,0,0}^{2}} \\
   {{r}_{1,2}}=& 72 {{a}_{3,1,0}} {{a}_{4,0,0}}+36 {{a}_{2,2,0}} {{a}_{3,1,0}}+18 {{a}_{2,1,1}} {{a}_{3,0,1}}+36 {{a}_{1,3,0}} {{a}_{2,2,0}}+ \\& 12 {{a}_{1,2,1}} {{a}_{2,1,1}}+12 {{a}_{1,1,2}} {{a}_{2,0,2}}+72 {{a}_{0,4,0}} {{a}_{1,3,0}}+18 {{a}_{0,3,1}} {{a}_{1,2,1}}+\\& 12 {{a}_{0,2,2}} {{a}_{1,1,2}}+18 {{a}_{0,1,3}} {{a}_{1,0,3}}\\ 
   {{r}_{1,3}}=& 72 {{a}_{3,0,1}} {{a}_{4,0,0}}+18 {{a}_{2,1,1}} {{a}_{3,1,0}}+36 {{a}_{2,0,2}} {{a}_{3,0,1}}+12 {{a}_{1,2,1}} {{a}_{2,2,0}}+\\& 12 {{a}_{1,1,2}} {{a}_{2,1,1}}+36 {{a}_{1,0,3}} {{a}_{2,0,2}}+18 {{a}_{0,3,1}} {{a}_{1,3,0}}+12 {{a}_{0,2,2}} {{a}_{1,2,1}}+\\& 18 {{a}_{0,1,3}} {{a}_{1,1,2}}+72 {{a}_{0,0,4}} {{a}_{1,0,3}}\\  
   {{r}_{2,2}} =&-96 {{a}_{4,0,0}^{2}}-6 {{a}_{3,1,0}^{2}}-24 {{a}_{3,0,1}^{2}}+8 {{a}_{2,2,0}^{2}}-2 {{a}_{2,1,1}^{2}}-16 {{a}_{2,0,2}^{2}}+\\& 30 {{a}_{1,3,0}^{2}}+4 {{a}_{1,2,1}^{2}}-2 {{a}_{1,1,2}^{2}}-24 {{a}_{1,0,3}^{2}}+192 {{a}_{0,4,0}^{2}}+30 {{a}_{0,3,1}^{2}}+\\& 8 {{a}_{0,2,2}^{2}}-6 {{a}_{0,1,3}^{2}}-96 {{a}_{0,0,4}^{2}} \\ 
   {{r}_{2,3}}= &18 {{a}_{3,0,1}} {{a}_{3,1,0}}+12 {{a}_{2,1,1}} {{a}_{2,2,0}}+12 {{a}_{2,0,2}} {{a}_{2,1,1}}+18 {{a}_{1,2,1}} {{a}_{1,3,0}}+\\& 12 {{a}_{1,1,2}} {{a}_{1,2,1}}+18 {{a}_{1,0,3}} {{a}_{1,1,2}}+72 {{a}_{0,3,1}} {{a}_{0,4,0}}+36 {{a}_{0,2,2}} {{a}_{0,3,1}}+\\& 36 {{a}_{0,1,3}} {{a}_{0,2,2}}+72 {{a}_{0,0,4}} {{a}_{0,1,3}}\\ 
  {{r}_{3,3}}= &-96 {{a}_{4,0,0}^{2}}-24 {{a}_{3,1,0}^{2}}-6 {{a}_{3,0,1}^{2}}-16 {{a}_{2,2,0}^{2}}-2 {{a}_{2,1,1}^{2}}+8 {{a}_{2,0,2}^{2}}\\& -24 {{a}_{1,3,0}^{2}}-2 {{a}_{1,2,1}^{2}}+4 {{a}_{1,1,2}^{2}}+30 {{a}_{1,0,3}^{2}}-96 {{a}_{0,4,0}^{2}}-6 {{a}_{0,3,1}^{2}}+\\& 8 {{a}_{0,2,2}^{2}}+30 {{a}_{0,1,3}^{2}}+192 {{a}_{0,0,4}^{2}}
\end{align*}

Functions and notations defined as above, we have 
${\mathscr{T}_2}$ is

\begin{align*}
& \{\, {{x}^{3}}y+{{x}^{4}},{{x}^{2}}\,
{{y}^{2}}+{{x}^{4}},{{x}^{2}}\,
{{y}^{2}}+{{x}^{3}}y,x\,{{y}^{3}}+{{x}^{4}},x\,{{y}^{3}}+{{x}^{3}}y,
{{y}^{4}}+{{x}^{4}},\\ &  
{{x}^{3}}z+{{x}^{3}}y, 
{{x}^{3}}z+{{x}^{2}}\,{{y}^{2}},
{{x}^{3}}z+x\,{{y}^{3}},
{{x}^{3}}z+{{y}^{4}},
{{x}^{2}}yz+{{x}^{4}}, \\ & 
{{x}^{2}}yz+{{x}^{3}}y,
{{x}^{2}}yz+{{x}^{2}}\,{{y}^{2}},{{x}^{2}}yz+x\,{{y}^{3}},
{{x}^{2}}yz+{{y}^{4}},
x\,{{y}^{2}}z+{{x}^{3}}z,\\ & 
x\,{{y}^{2}}z+{{x}^{2}}yz,
{{y}^{3}}z+{{x}^{3}}z,
{{x}^{2}}\,{{z}^{2}}+{{x}^{2}}\,{{y}^{2}},
\\ & {{x}^{2}}\,{{z}^{2}}+x\,{{y}^{3}}
,{{x}^{2}}\,{{z}^{2}}+{{y}^{4}},
{{x}^{2}}\,{{z}^{2}}+x\,{{y}^{2}}z\}\end{align*}
and ${\mathscr{T}_3}$ is
\begin{align*}
&\{\, {{x}^{2}}\,{{y}^{2}}+{{x}^{3}}y+{{x}^{4}},x\,{{y}^{3}}+{{x}^{3}}y+{{x}^{4}},x\,{{y}^{3}}+{{x}^{2}}\,{{y}^{2}}+{{x}^{4}},\\
&  x\,{{y}^{3}}+{{x}^{2}}\,{{y}^{2}}+{{x}^{3}}y, {{y}^{4}}+{{x}^{3}}y+{{x}^{4}},{{y}^{4}}+{{x}^{2}}\,{{y}^{2}}+{{x}^{4}}, \\ &{{x}^{3}}z+{{x}^{3}}y+{{x}^{4}},{{x}^{3}}z+{{x}^{2}}\,{{y}^{2}}+{{x}^{4}}, {{x}^{3}}z+{{x}^{2}}\,{{y}^{2}}+{{x}^{3}}y, \\ &{{x}^{3}}z+x\,{{y}^{3}}+{{x}^{4}},{{x}^{3}}z+x\,{{y}^{3}}+{{x}^{3}}y, {{x}^{3}}z+x\,{{y}^{3}}+{{x}^{2}}\,{{y}^{2}}, \\ &{{x}^{3}}z+{{y}^{4}}+{{x}^{4}},{{x}^{3}}z+{{y}^{4}}+{{x}^{3}}y,{{x}^{3}}z+{{y}^{4}}+{{x}^{2}}\,{{y}^{2}}, \\ & {{x}^{3}}z+{{y}^{4}}+x\,{{y}^{3}},{{x}^{2}}yz+{{x}^{3}}y+{{x}^{4}},{{x}^{2}}yz+{{x}^{2}}\,{{y}^{2}}+{{x}^{4}}, \\ &{{x}^{2}}yz+{{x}^{2}}\,{{y}^{2}}+{{x}^{3}}y, {{x}^{2}}yz+x\,{{y}^{3}}+{{x}^{4}},{{x}^{2}}yz+x\,{{y}^{3}}+{{x}^{3}}y, \\ &{{x}^{2}}yz+x\,{{y}^{3}}+{{x}^{2}}\,{{y}^{2}}, {{x}^{2}}yz+{{y}^{4}}+{{x}^{4}},{{x}^{2}}yz+{{y}^{4}}+{{x}^{3}}y, \\ &{{x}^{2}}yz+{{y}^{4}}+{{x}^{2}}\,{{y}^{2}},{{x}^{2}}yz+{{y}^{4}}+x\,{{y}^{3}}, {{x}^{2}}yz+{{x}^{3}}z+{{x}^{3}}y, \\ &{{x}^{2}}yz+{{x}^{3}}z+{{x}^{2}}\,{{y}^{2}},{{x}^{2}}yz+{{x}^{3}}z+x\,{{y}^{3}},{{x}^{2}}yz+{{x}^{3}}z+{{y}^{4}},  \\  &x\, {{y}^{2}}z+{{x}^{3}}z+{{x}^{4}},x\,{{y}^{2}}z+{{x}^{3}}z+{{x}^{3}}y,x\,{{y}^{2}}z+{{x}^{3}}z+{{x}^{2}}\,{{y}^{2}}, \\ &x\,{{y}^{2}}z+{{x}^{3}}z+x\,{{y}^{3}},x\,{{y}^{2}}z+{{x}^{3}}z+{{y}^{4}},
 x\,{{y}^{2}}z+{{x}^{2}}yz+{{x}^{4}}, \\ &x\,{{y}^{2}}z+{{x}^{2}}yz+{{x}^{3}}y,x\,{{y}^{2}}z+{{x}^{2}}yz+{{x}^{2}}\,{{y}^{2}},x\,{{y}^{2}}z+{{x}^{2}}yz+{{x}^{3}}z, \\ & {{y}^{3}}z+{{x}^{3}}z+{{x}^{4}},{{y}^{3}}z+{{x}^{3}}z+{{x}^{3}}y,{{y}^{3}}z+{{x}^{3}}z+{{x}^{2}}\,{{y}^{2}}, \\ &{{y}^{3}}z+{{x}^{2}}yz+{{x}^{3}}z,{{x}^{2}}\,{{z}^{2}}+{{x}^{2}}\,{{y}^{2}}+{{x}^{4}},  {{x}^{2}}\,{{z}^{2}}+{{x}^{2}}\,{{y}^{2}}+{{x}^{3}}y, \\ \end{align*} \begin{align*} &{{x}^{2}}\,{{z}^{2}}+x\,{{y}^{3}}+{{x}^{4}},{{x}^{2}}\,{{z}^{2}}+x\,{{y}^{3}}+{{x}^{3}}y,  {{x}^{2}}\,{{z}^{2}}+x\,{{y}^{3}}+{{x}^{2}}\,{{y}^{2}}, \\ &{{x}^{2}}\,{{z}^{2}}+{{y}^{4}}+{{x}^{4}},{{x}^{2}}\,{{z}^{2}}+{{y}^{4}}+{{x}^{3}}y, {{x}^{2}}\,{{z}^{2}}+{{y}^{4}}+{{x}^{2}}\,{{y}^{2}}, \\ &{{x}^{2}}\,{{z}^{2}}+{{y}^{4}}+x\,{{y}^{3}},{{x}^{2}}\,{{z}^{2}}+{{x}^{3}}z+x\,{{y}^{3}},  {{x}^{2}}\,{{z}^{2}}+{{x}^{3}}z+{{y}^{4}}, \\ &{{x}^{2}}\,{{z}^{2}}+{{x}^{2}}yz+{{x}^{2}}\,{{y}^{2}},{{x}^{2}}\,{{z}^{2}}+{{x}^{2}}yz+x\,{{y}^{3}},  {{x}^{2}}\,{{z}^{2}}+{{x}^{2}}yz+{{y}^{4}}, \\ &{{x}^{2}}\,{{z}^{2}}+x\,{{y}^{2}}z+{{x}^{4}},{{x}^{2}}\,{{z}^{2}}+x\,{{y}^{2}}z+{{x}^{3}}y, {{x}^{2}}\,{{z}^{2}}+x\,{{y}^{2}}z+{{x}^{2}}\,{{y}^{2}}, \\ &{{x}^{2}}\,{{z}^{2}}+x\,{{y}^{2}}z+x\,{{y}^{3}},{{x}^{2}}\,{{z}^{2}}+x\,{{y}^{2}}z+{{y}^{4}},  {{x}^{2}}\,{{z}^{2}}+x\,{{y}^{2}}z+{{x}^{3}}z, \\ &{{x}^{2}}\,{{z}^{2}}+x\,{{y}^{2}}z+{{x}^{2}}yz,{{x}^{2}}\,{{z}^{2}}+{{y}^{3}}z+{{x}^{4}},  {{x}^{2}}\,{{z}^{2}}+{{y}^{3}}z+{{x}^{3}}y, \\ &{{x}^{2}}\,{{z}^{2}}+{{y}^{3}}z+{{x}^{2}}\,{{y}^{2}},{{x}^{2}}\,{{z}^{2}}+{{y}^{3}}z+x\,{{y}^{3}},  {{x}^{2}}\,{{z}^{2}}+{{y}^{3}}z+{{x}^{3}}z, \\ &{{x}^{2}}\,{{z}^{2}}+{{y}^{3}}z+{{x}^{2}}yz,xy\,{{z}^{2}}+x\,{{y}^{3}}+{{x}^{4}},  xy\,{{z}^{2}}+x\,{{y}^{3}}+{{x}^{3}}y, \\ &xy\,{{z}^{2}}+{{y}^{4}}+{{x}^{4}},xy\,{{z}^{2}}+{{x}^{3}}z+x\,{{y}^{3}},  xy\,{{z}^{2}}+{{x}^{3}}z+{{y}^{4}}, \\ &xy\,{{z}^{2}}+{{x}^{2}}yz+x\,{{y}^{3}},xy\,{{z}^{2}}+{{x}^{2}}yz+{{y}^{4}},xy\,{{z}^{2}}+x\,{{y}^{2}}z+{{x}^{2}}yz, \\ &  xy\,{{z}^{2}}+{{y}^{3}}z+{{x}^{3}}z,{{y}^{2}}\,{{z}^{2}}+{{x}^{2}}\,{{z}^{2}}+{{x}^{2}}\,{{y}^{2}},x\,{{z}^{3}}+x\,{{y}^{3}}+{{x}^{4}}, \\ &x\,{{z}^{3}}+x\,{{y}^{3}}+{{x}^{3}}y,  x\,{{z}^{3}}+{{y}^{4}}+{{x}^{4}},x\,{{z}^{3}}+{{y}^{4}}+{{x}^{3}}y, \\ &x\,{{z}^{3}}+{{x}^{3}}z+{{y}^{4}},x\,{{z}^{3}}+{{y}^{3}}z+{{x}^{3}}y,{{z}^{4}}+{{y}^{4}}+{{x}^{4}}\, \}
\end{align*}

Then the  polynomials with diagonal moment matrices are
\begin{align*}& 
{{b}_{1}}\,{{x}^{3}}z+{{b}_{2}}\,{{y}^{4}}+{{x}^{2}}\,{{y}^{2}},{{b}_{1}}\,{{x}^{2}}yz+{{b}_{2}}x\,{{y}^{3}}+{{x}^{4}}, \\ & {{b}_{1}}\,{{x}^{2}}yz+{{b}_{2}}\,{{y}^{4}}+{{x}^{4}},{{b}_{1}}x\,{{y}^{2}}z+{{b}_{2}}\,{{x}^{3}}z+{{y}^{4}},\\
& {{b}_{1}}\,{{y}^{3}}z+{{b}_{2}}\,{{x}^{3}}z+{{x}^{2}}\,{{y}^{2}}, {{b}_{1}}\,{{x}^{2}}\,{{z}^{2}}+{{b}_{2}}\,{{x}^{2}}\,{{y}^{2}}+{{x}^{4}},\\
& {{b}_{1}}\,{{x}^{2}}\,{{z}^{2}}+{{b}_{2}}x\,{{y}^{3}}+{{x}^{4}},{{b}_{1}}\,{{x}^{2}}\,{{z}^{2}}+{{b}_{2}}x\,{{y}^{3}}+{{x}^{3}}y, \\ & {{b}_{1}}\,{{x}^{2}}\,{{z}^{2}}+{{b}_{2}}\,{{y}^{4}}+{{x}^{4}},{{b}_{1}}\,{{x}^{2}}\,{{z}^{2}}+{{b}_{2}}\,{{y}^{4}}+{{x}^{3}}y,\\ 
& {{b}_{1}}\,{{x}^{2}}\,{{z}^{2}}+{{b}_{2}}\,{{y}^{4}}+{{x}^{2}}\,{{y}^{2}}, {{b}_{1}}\,{{x}^{2}}\,{{z}^{2}}+{{b}_{2}}x\,{{y}^{2}}z+{{x}^{4}},\\
& {{b}_{1}}\,{{x}^{2}}\,{{z}^{2}}+{{b}_{2}}x\,{{y}^{2}}z+{{x}^{3}}y, {{b}_{1}}\,{{x}^{2}}\,{{z}^{2}}+{{b}_{2}}x\,{{y}^{2}}z+{{y}^{4}},\\
& {{b}_{1}}\,{{x}^{2}}\,{{z}^{2}}+{{b}_{2}}\,{{y}^{3}}z+{{x}^{4}},{{b}_{1}}\,{{x}^{2}}\,{{z}^{2}}+{{b}_{2}}\,{{y}^{3}}z+{{x}^{3}}y,\\ & {{b}_{1}}\,{{x}^{2}}\,{{z}^{2}}+{{b}_{2}}\,{{y}^{3}}z+{{x}^{2}}\,{{y}^{2}},{{b}_{1}}xy\,{{z}^{2}}+{{b}_{2}}x\,{{y}^{3}}+{{x}^{4}},\\ & {{b}_{1}}xy\,{{z}^{2}}+{{b}_{2}}x\,{{y}^{3}}+{{x}^{3}}y,{{b}_{1}}xy\,{{z}^{2}}+{{b}_{2}}\,{{y}^{4}}+{{x}^{4}},\\ & {{b}_{1}}xy\,{{z}^{2}}+{{b}_{2}}\,{{x}^{3}}z+x\,{{y}^{3}},{{b}_{1}}xy\,{{z}^{2}}+{{b}_{2}}\,{{x}^{3}}z+{{y}^{4}},\\ & {{b}_{1}}xy\,{{z}^{2}}+{{b}_{2}}\,{{y}^{3}}z+{{x}^{3}}z,{{b}_{1}}\,{{y}^{2}}\,{{z}^{2}}+{{b}_{2}}\,{{x}^{2}}\,{{z}^{2}}+{{x}^{2}}\,{{y}^{2}},\\ & {{b}_{1}}x\,{{z}^{3}}+{{b}_{2}}x\,{{y}^{3}}+{{x}^{4}},{{b}_{1}}x\,{{z}^{3}}+{{b}_{2}}x\,{{y}^{3}}+{{x}^{3}}y,\\
& {{b}_{1}}x\,{{z}^{3}}+{{b}_{2}}\,{{y}^{4}}+{{x}^{4}},  {{b}_{1}}x\,{{z}^{3}}+{{b}_{2}}\,{{y}^{4}}+{{x}^{3}}y,{{b}_{1}}x\,{{z}^{3}}+{{b}_{2}}\,{{x}^{3}}z+{{y}^{4}},\\ & {{b}_{1}}x\,{{z}^{3}}+{{b}_{2}}\,{{y}^{3}}z+{{x}^{3}}y,{{b}_{1}}\,{{z}^{4}}+{{b}_{2}}\,{{y}^{4}}+{{x}^{4}}
\end{align*}

Here we explain more details on how to find the critical points of the square length of the moment map.
We use 
\begin{verbatim}
block( 
    C3:sublist(T(3),
       lambda([x],diagmatrixp(Moment(x),n))),
    CRT3:create_list(algsys(
         create_list(i[j]=0,j,1,B),[b[1],b[2]]),i,
      create_list(ratsimp(Grad(i)),i,C3)),
    RCRT3:create_list(sublist(i,IsReal),i,CRT3)
  );
\end{verbatim}
to find solutions for the $b_i$'s such that the polynomial has the gradient of the length of its moment matrix vanishes. Here $\tt C3$ chooses elements in $\tt C3$ whose the moment matrices are diagonal. $\tt CRT3$ solves the problem of being a critical point of the square length of the moment matrix in $\tt C3$. Since we are solving this problem in the reals ${\mathbb R}$, all complex solutions should not be involved, and we use $\tt RCRT3$ to select real solutions. 
Then we use
\begin{verbatim}
CC3:[]\$;
CRTV3:for i:1 thru length(RCRT3) do 
        for j in RCRT3[i] do  
           push(subst(j,C3[i]),CC3),CC3;
CCC3:sublist(CC3,fff);
\end{verbatim}
to substitute the solutions back into the polynomials. At the beginning $\tt CC3$ is an empty list, and for each time we push the substitution into $\tt CC3$. But the function ``{\tt algsys}" sometimes gives incorrect answers for a system of equations. So we have to check our answer by
\begin{verbatim}
CCCC3:sublist(CCC3,lambda([x], is(
      Grad(x)=[0,0,0,0,0,0,0,0,0,0,0,0,0,0,0])));
\end{verbatim}

Here are the solutions of the critical points for $l(f)\leq 3$ up to obvious isomorphism (a scalar product on one of the coordinates). We draw the affine cones of them in Mathematica.

These pictures help us to see the singularities of the curves. Nevertheless, those empty graphs remind us that some of these curves do not have real solutions, because the only solutions of in ${\mathbb R}^3$ is the point $(0,0,0)$ which is not contained in the affine cone of the projective plane, and we should delete them from the real critical points.

\begin{proposition}The critical points for $l(f)\leq 3$ up to obvious isomorphism (a scalar product on one of the coordinates) are
\begin{equation}
\label{critquartic}\begin{split}
& {{x}^{3}} z+x {{y}^{3}},\frac{{{x}^{3}} z}{\sqrt{3}}+{{x}^{2}} {{y}^{2}},3\cdot{{2}^{\frac{3}{2}}} x {{y}^{2}} z+{{x}^{4}},\sqrt{15} x {{y}^{2}} z+{{x}^{3}} y,\\
& {{y}^{3}} z+{{x}^{3}} y,  {{y}^{3}} z+{{x}^{3}} z,\frac{{{y}^{3}} z}{\sqrt{3}}+{{x}^{2}} {{y}^{2}},\\
& \sqrt{3} {{x}^{2}} {{z}^{2}}+{{x}^{3}} y,\sqrt{3} {{x}^{2}} {{z}^{2}}+\frac{2 {{y}^{3}} z}{\sqrt{3}}+{{x}^{2}} {{y}^{2}},  \\
& 2 \sqrt{6} {{x}^{2}} {{z}^{2}}+4 {{y}^{3}} z+{{x}^{4}},\frac{\sqrt{21} {{x}^{2}} {{z}^{2}}}{2}+\sqrt{3} {{y}^{3}} z+{{x}^{3}} y,\\
& 3 x y {{z}^{2}}+{{x}^{3}} y,  3 x y {{z}^{2}}+x {{y}^{3}},3\cdot {{2}^{\frac{3}{2}}} x y {{z}^{2}}+{{x}^{4}}, \\
& \\
& 3\cdot {{2}^{\frac{3}{2}}} x y {{z}^{2}}+{{2}^{\frac{3}{2}}} x {{y}^{3}}+{{x}^{4}},3\cdot {{2}^{\frac{3}{2}}} x y {{z}^{2}}+{{y}^{4}},\\
 & {{2}^{\frac{5}{2}}} x y {{z}^{2}}+\frac{4 {{x}^{3}} z}{\sqrt{3}}+{{y}^{4}},\sqrt{3} x y {{z}^{2}}+{{y}^{3}} z+{{x}^{3}} z,\\
& 2 \sqrt{3} x y {{z}^{2}}+x {{y}^{3}}+{{x}^{3}} y,
 4 \sqrt{3} x y {{z}^{2}}+{{y}^{4}}+{{x}^{4}}, \\
& \sqrt{7} x y {{z}^{2}}+\frac{\sqrt{2} {{x}^{3}} z}{\sqrt{3}}+x {{y}^{3}},\sqrt{15} x y {{z}^{2}}+{{x}^{3}} z,\\
& x {{z}^{3}}+{{x}^{3}} y,x {{z}^{3}}+{{y}^{3}} z+{{x}^{3}} y,\\
& \\
& 2 x {{z}^{3}}+2 {{x}^{3}} z+{{y}^{4}},  4 x {{z}^{3}}+4 x {{y}^{3}}+{{x}^{4}}\end{split}
\end{equation}

\end{proposition}

\end{document}